\documentclass[12pt,a4paper]{article}
\usepackage{mathrsfs}
\usepackage{amssymb}

\usepackage{amsmath}

\newtheorem{theorem}{Theorem}[section]
\newtheorem{definition}[theorem]{Definition}
\newtheorem{lemma}[theorem]{Lemma}

\begin{document}

\title{Gr\"{o}bner-Shirshov bases for  extensions of
algebras\footnote {The research is supported by the National Natural
Science Foundation of China (Grant No.10771077) and the Natural
Science Foundation of Guangdong Province (Grant No.06025062).} }
\author{
Yuqun Chen  \\
\\
{\small \ School of Mathematical Sciences}\\
{\small \ South China Normal University}\\
{\small \ Guangzhou 510631}\\
{\small \ P. R. China}\\
{\small \ yqchen@scnu.edu.cn} }
\date{}

\maketitle \noindent{\bf Abstract.}
  An algebra $\cal{R}$ is called an extension of the
algebra $M$ by $B$ if $M^2=0$, $M$ is an ideal of $\cal{R}$ and
$\cal{R}$$/M\cong B$ as algebras. In this paper, by using the
Gr\"{o}bner-Shirshov bases,  we characterize completely the
extensions of  $M$ by $B$. An algorithm to find the conditions of an
algebra $A$ to be an extension of $M$ by $B$ is obtained.

\section{Preliminaries}

Let $k$ be a field, $X$ a set, $X^{*} $ the monoid of all words on
$X$ and $X^{+}$ the free semigroup of nonempty words on $X$. We
denote $k\langle X_+ \rangle$ the $k$-span of all nonempty words
in $X$. As we know, $k\langle X_+ \rangle$ is a free associative
algebra without identity and $k\langle X\rangle$ a free
associative algebra with identity.  It is clear that, for every
algebra $\mathscr{A}$ (not necessarily with 1), we have
$\mathscr{A} \cong k\langle X_+ \rangle /I$ for some $X$ and ideal
$I$ of $ k \langle X_+ \rangle$. For a word $w\in X^*$, we denote
the length of $w$ by $deg(w)$. Let $X^*$ be a well ordered set.
Let $f\in k\langle X\rangle$ with the leading word $\bar{f}$. We
say that $f$ is monic if $\bar{f}$ has coefficient 1.

\begin{definition} (\cite{Sh}, see also \cite{b72}, \cite{b76}) \
Let $f$ and $g$ be two monic polynomials in \textmd{k}$\langle
X\rangle$ and $<$ a well order on $X^*$. Then, there are two kinds
of compositions:

$(1)$ If \ $w$ is a word such that $w=\bar{f}b=a\bar{g}$ for some
$a,b\in X^*$ with deg$(\bar{f})+$deg$(\bar{g})>$deg$(w)$, then the
polynomial
 $(f,g)_w=fb-ag$ is called the intersection composition of $f$ and
$g$ with respect to $w$.

$(2)$ If  $w=\bar{f}=a\bar{g}b$ for some $a,b\in X^*$, then the
polynomial $(f,g)_w=f - agb$ is called the inclusion composition
of $f$ and $g$ with respect to $w$.

\end{definition}

\begin{definition}(\cite{b72}, \cite{b76}, cf. \cite{Sh})
Let $S\subseteq$ $\textmd{k}\langle X\rangle$ with each $s\in S$
monic. Then the composition $(f,g)_w$ is called trivial modulo $S$
if $(f,g)_w=\sum\alpha_i a_i s_i b_i$, where each $\alpha_i\in k$,
$a_i,b_i\in X^{*}$ and $\overline{a_i s_i b_i}<w$. If this is the
case, then we write
$$
(f,g)_w\equiv0\quad mod(S,w)
$$
In general, for $p,q\in k\langle X\rangle$, we write
$$
p\equiv q\quad mod(S,w)
$$
which means that $p-q=\sum\alpha_i a_i s_i b_i $, where each
$\alpha_i\in k,a_i,b_i\in X^{*}$ and $\overline{a_i s_i b_i}<w$.
\end{definition}

\begin{definition} (\cite{b72}, \cite{b76}, cf. \cite{Sh}) \
We call the set $S$ with respect to the well order a
Gr\"{o}bner-Shirshov set (basis) in $k\langle X\rangle$ if any
composition of polynomials in $S$ is trivial relative to $S$.
\end{definition}

A well order $>$ on $X^*$ is monomial if it is compatible with the
multiplication of words, that is, for $u, v\in X^*$, we have
$$
u > v \Rightarrow w_{1}uw_{2} > w_{1}vw_{2},  \ for \  all \
 w_{1}, \ w_{2}\in  X^*.
$$
A standard example of monomial order on $X^*$ is the deg-lex order
to compare two words first by degree and then lexicographically,
where $X$ is a linearly ordered set.

The following lemma was proved by Shirshov \cite{Sh} for the free
Lie algebras (with deg-lex ordering) in 1962 (see also Bokut
\cite{b72}). In 1976, Bokut \cite{b76} specialized the approach of
Shirshov to associative algebras (see also Bergman \cite{b}). For
commutative polynomials, this lemma is known as the Buchberger's
Theorem in \cite{bu65} and \cite{bu70}.

\begin{lemma}\label{l1}
(Composition-Diamond Lemma) \ Let $k$ be a field, $A=k \langle
X|S\rangle=K\langle X\rangle/Id(S)$ and $<$ a monomial order on
$X^*$, where $Id(S)$ is the ideal of $k \langle X\rangle$
generated by $S$. Then the following statements are equivalent:
\begin{enumerate}
\item[(i)] $S $ is a Gr\"{o}bner-Shirshov basis.
\item[(ii)] $f\in
Id(S)\Rightarrow \bar{f}=a\bar{s}b$ for some $s\in S$ and $a,b\in
X^*$.
\item[(iii)] $Irr(S) = \{ u \in X^* |  u \neq a\bar{s}b
,s\in S,a ,b \in X^*\}$ is a basis of the algebra $A=k\langle X |
S \rangle$.
\end{enumerate}
\end{lemma}

The following lemma comes from \cite{cq} (\cite{cq}, Lemma 4.2)
which is essentially the same as Lemma \ref{l1}.

\begin{lemma}\label{l2}
{\em (Composition-Diamond Lemma)} \ Let $k$ be a field, $S\subseteq
K\langle X_+\rangle$ be  monic, $A=k \langle X_+|S\rangle=K\langle
X_+\rangle/Id(S)$ and $<$ a monomial order on $X^+$, where $Id(S)$
is the ideal of $k \langle X_+\rangle$ generated by $S$. Then the
following statements are equivalent:
\begin{enumerate}
\item[(i)] $S $ is a Gr\"{o}bner-Shirshov basis.
\item[(ii)] $f\in
Id(S)\Rightarrow \bar{f}=a\bar{s}b$ for some $s\in S$ and $a,b\in
X^*$.
\item[(iii)] $Irr(S) = \{ u \in X^+ |  u \neq a\bar{s}b
,s\in S,a ,b \in X^*\}$ is a basis of the algebra $A=k\langle X_+
| S \rangle$.
\end{enumerate}
\end{lemma}

\noindent{\bf Remark:} Suppose that $S\subseteq K\langle X_+\rangle$
is a Gr\"{o}bner-Shirshov  basis. In the (iii) of Lemma \ref{l2},
$1\notin Irr(S)$ but not the case in the Lemma \ref{l1}.

\ \

For convenience, we identify a relation $u=f_u$ of an algebra
presented by generators and relations with the polynomial $u-f_u$
in the corresponding free algebra.

The concepts of the extension of algebras was invented by Hochschild
\cite{h} (also see \cite{ce} and \cite{m}).

\begin{definition} Let $k$ be a field, $M,B, \cal{R}$ $k$-algebras (not necessarily with 1).
Then $\cal{R}$ is called an extension of $M$ by $B$ if  $M^2=0$,
$M$ is an ideal of $\cal{R}$ and $\cal{R}$$/M\cong B$ as algebras.
\end{definition}

In \cite{h}, such an extension is called a singular extension.

\section{Characterizations of extensions of
algebras}

Let $M,B$ be $k$-algebras, $M^2=0$ and $M$ a $B$-bimodule.

By a factor set $\{(b,b')|b,b'\in B\}$  of $B$ in $M$ we mean that
$\{(b,b')|b,b'\in B\}$ is a subset of $M$ such that the function
$(b,b')$ is $k$-bilinear.

Let $I, \ J$ be linearly ordered sets, $\{b_i|i\in I\}, \
\{m_j|j\in J\} \ k$-bases of $B$ and $M$, respectively and
$\{(b,b')|b,b'\in B\}$ a factor set of $B$ in $M$. Denote
$$
A=E_k(M,B,(b_p,b_q))=k\langle (\{m_{_j}\},\{b_i\})_+| S \rangle
$$
where $S=\{b_pb_q=[b_pb_q]+(b_p,b_q),b_pm_j=[b_pm_j],
m_jb_p=[m_jb_p],m_jm_l=0,p,q\in I,j,l\in J\}$ and for example,
$[b_pb_q]=\sum_{i\in I}\alpha_{pq}^ib_i, \
 \alpha_{pq}^i\in k,$ is the product in $B$.

We order the set $(\{m_j\}_J\cup\{b_i\}_I)^+$ by the deg-lex order.

Equipping the above concepts, we have the following theorems which
give characterizations of  extensions of algebras.

\begin{theorem}\label{st1}
Let $M,B$ be $k$-algebras, $M^2=0, \ M$ a $B$-bimodule,
$\{b_i|i\in I\}, \ \{m_j|j \in J\}$  $k$-bases of $B$ and $M$,
respectively and $\{(b,b')|b,b'\in B\}$ a factor set of $B$ in
$M$. Let $A=E_k(M,B,(b_p,b_q))=k\langle (\{m_{_j}\},\{b_i\})_+| S
\rangle$, where $S=\{b_pb_q=[b_pb_q]+(b_p,b_q),b_pm_j=[b_pm_j],
m_jb_p=[m_jb_p],m_jm_l=0,p,q\in I,j,l\in J\}$. Then, $S$ is a
Gr\"{o}bner-Shirshov bases for $A$ if and only if the factor set
satisfies the following condition in $M$: for any $p,q,r\in I$
\begin{eqnarray}\label{e1}
 b_p(b_q,b_r)-(b_pb_q,b_r)+(b_p,b_qb_r)-(b_p,b_q)b_r=0
\end{eqnarray}
i.e., the function $(b,b')$ is a cocycle.

 Moreover, if this
is the case, $A=E_k(M,B,(b_p,b_q))= k\langle (\{m_{_j}\},\{b_i\})_+|
S \rangle$ is an extension of $M$ by $B$ in a natural way.
\end{theorem}
\noindent {\bf Proof.} \ Suppose that (\ref{e1}) holds. Then the
possible compositions in $S$ are related to the following
ambiguities:
$$
w=b_pb_qb_r, \ m_jb_pb_q, \ b_pb_qm_j, \ m_jm_lm_n, \ m_jm_lb_p, \
b_pm_jm_l, \ b_pm_jb_q, \ m_jb_pm_l.
$$

For $w=m_jm_lm_n, \ m_jm_lb_p, \ b_pm_jm_l, \ b_pm_jb_q, \
m_jb_pm_l$, by noting that $[[m_jm_l]m_n]=[m_j[m_lm_n]]=0, \
[[m_jm_l]b_p]=[m_j[m_lb_p]]=0, \ [[b_pm_j]m_l]=[b_p[m_jm_l]]=0, \
[[m_jb_p]m_l]=[m_j[b_pm_l]]=0$ and $[[b_pm_j]b_q]=[b_p[m_jb_q]]$ for
any $j,l,n\in J, \ p,q\in I$, the correspondent compositions are
trivial modulo $S$.

For $w=m_jb_pb_q$, since  $[[m_jb_p]b_q]=[m_j[b_pb_q]]$ and
$m_j(b_p,b_q)=0$, we have
\begin{eqnarray*}
&&(m_jb_p-[m_jb_p],b_pb_q-[b_pb_q]-(b_p,b_q))_w=
-[m_jb_p]b_q+m_j([b_pb_q]+(b_p,b_q))\\
&\equiv& -[[m_jb_p]b_q]+m_j[b_pb_q]+m_j(b_p,b_q) \equiv
-[[m_jb_p]b_q]+[m_j[b_pb_q]]\equiv 0 \  \ mod(S,w)
\end{eqnarray*}

Similarly, for $w=b_pb_qm_j$, the correspondent composition is
trivial modulo $S$.

For $w=b_pb_qb_r$, by noting that $[[b_pb_q]b_r]=[b_p[b_qb_r]]$, we have
\begin{eqnarray*}
&&(b_pb_q-[b_pb_q]-(b_p,b_q),b_qb_r-[b_qb_r]-(b_q,b_r))_w\\
&=&
-([b_pb_q]+(b_p,b_q))b_r+b_p([b_qb_r]+(b_q,b_r))\\
&\equiv& -[[b_pb_q]b_r]-(b_pb_q,b_r)-(b_p,b_q)b_r+([b_p[b_qb_r]]
+(b_p,b_qb_r)+b_p(b_q,b_r))\\
&\equiv&0 \  \ mod(S,w)
\end{eqnarray*}

Thus, all compositions in $S$ are trivial modulo $S$.

Conversely, if $S$ is a Gr\"{o}bner-Shirshov bases for $A$, then,
by noting that for $w=b_pb_qb_r$, the compositions are trivial
modulo $S$, we can easily check that the condition (\ref{e1})
holds in $\overline{M}=(M+Id(S))/Id(S)$. Then, by using Lemma
\ref{l2}, the algebra $\overline{M}$ is the same as $M$.

Moreover, assume that $S$ is a Gr\"{o}bner-Shirshov bases for $A$.
Then, by Lemma \ref{l2}, each element $r\in A$ can be uniquely
written as $r=m+b$, where $m\in M, \ b\in B$. Hence, $A=M\oplus B$
as $k$-modules with the following multiplication: for any $m,m'\in
M, \ b,b'\in B$,
$$
(m+b)\cdot (m'+b')=mb'+bm'+(b,b')+bb'.
$$
From this it follows that $M$ is an ideal of $A$ and $A/M\cong B$ as
algebras. \ \ $\square$

\ \

\begin{theorem}\label{st2}
Let $M,B, \cal{R}$ be $k$-algebras with $M^2=0$. If $\cal{R}$ is an
extension of  $M$ by $B$ and  $\sigma: \ \cal{R}$$/M\rightarrow B,
r_b+M\mapsto b$ an algebra isomorphism, then $M$ is a $B$-bimodule
in a natural way: for any $b\in B, \ m\in M$,
\begin{eqnarray}\label{e5}
b\cdot m=r_bm, \ \ \ m\cdot b=mr_b
\end{eqnarray}
and there exists a factor set $\{(b,b')|b,b'\in B\}$ of $B$ in $M$
such that for any $b,b',b''\in B$,
$$
b(b',b'')-(bb',b'')+(b,b'b'')-(b,b')b''=0.
$$
Moreover, $\cal{R}$ $\cong A$ as algebras, where
$A=E_k(M,B,(b_p,b_q))= k\langle (\{m_{_j}\},\{b_i\})_+| S \rangle$,
$\{m_{_j}\}_J, \ \{b_i\}_I$ are any linear bases of $M$ and $B$
respectively and $S=\{b_pb_q=[b_pb_q]+(b_p,b_q),b_pm_j=[b_pm_j],
m_jb_p=[m_jb_p],m_jm_l=0,p,q\in I,j,l\in J\}$.
\end{theorem}
\noindent {\bf Proof.} \ Clearly, as $k$-modules, $\cal{R}$$=M\oplus
C$, where $C=\sum_{i\in I}kr_{b_i}$ with a basis $\{r_{b_i}|i\in
I\}$. Thus, for any $b,b'\in B$, we have $r_b+r_{b'}=r_{b+b'}$
because $\sigma$ is an algebra isomorphism. Since
$(r_b+M)(r_{b'}+M)=r_br_{b'}+M=r_{bb'}+M$, there exists a unique
$(b,b')\in M$ such that $r_br_{b'}=r_{bb'}+(b,b')$. Then, it is easy
to see that $M$ is a $B$-bimodule with the module operations
(\ref{e5}) and the function $(b,b')$ is $k$-bilinear. For example,
for any $b,b',b''\in B, m\in M, \ (bb')\cdot
m=r_{bb'}m=(r_br_{b'}-(b,b'))m=(r_br_{b'})m=r_b(r_{b'}m)=b\cdot(b'\cdot
m)$. Also, since
$r_b(r_{b'}+r_{b''})=r_br_{b'}+r_br_{b''}=r_br_{b'+b''}$, we have
$(b,b'+b'')=(b,b')+(b,b'')$.

Moreover, by noting that $(r_br_{b'})r_{b''}=r_b(r_{b'}r_{b''})$, we
know that the factor set $\{(b,b')|b,b'\in B\}$ satisfies
$$
b(b',b'')-(bb',b'')+(b,b'b'')-(b,b')b''=0.
$$

Now, we prove that $\cal{R}$ $\cong A$ as algebras. For any $m,m'\in
M, \ c,c'\in C$,
$$
(m+c)\cdot (m'+c')=mc'+cm'+(c,c')+cc'
$$
where $(c,c')=\sum_{p,q}\alpha_p\alpha_q'(b_p,b_q)$ if
$c=\sum_{p}\alpha_pr_{b_p}, \ c'=\sum_{q}\alpha_q'r_{b_q}$. Then, by
the proof of Theorem \ref{st1}, it is easy to see that $\tau: \
\cal{R}$$\rightarrow A$ by $\tau(m_j)=m_j, \ \tau(r_{b_i})=b_i$ is
an algebra isomorphism. \ \ $\square$

\ \

By Theorem \ref{st1} and Theorem \ref{st2}, we have the following
theorem.

\begin{theorem}\label{st6}
Let $M,B, \cal{R}$ be $k$-algebras with $M^2=0$. Then $\cal{R}$ is
an extension of  $M$ by $B$ if and only if $M$ is a $B$-bimodule and
there exists a factor set $\{(b,b')|b,b'\in B\}$ of $B$ in $M$ such
that for any $b,b',b''\in B$,
$$
b(b',b'')-(bb',b'')+(b,b'b'')-(b,b')b''=0
$$
and $\cal{R}$ $\cong A=E_k(M,B,(b_p,b_q))= k\langle
(\{m_{_j}\},\{b_i\})_+| S \rangle$, where $\{b_i\}_I$ and
$\{m_{_j}\}_J$ are any linear bases of $B$ and $M$ respectively,
$S=\{b_pb_q=[b_pb_q]+(b_p,b_q),b_pm_j=[b_pm_j],
m_jb_p=[m_jb_p],m_jm_l=0,p,q\in I,j,l\in J\}$.
\end{theorem}

\ \

Now, we consider the general case when the algebra $B$ is presented
by generators and relations.

Let $M,B$ be $k$-algebras, $M^2=0$, $M$ a $B$-bimodule,
$B=k\langle X_+|R \rangle$, where $R$ is a Gr\"{o}bner-Shirshov
bases for $B$ with the deg-lex order $<_B$ on $X^+$.

For convenience, we can assume that $R$ is a minimal
Gr\"{o}bner-Shirshov bases in a sense that the leading monomials
are not contained each other as subwords, in particular, they are
pairwise different.  Let $R=\{u-f_u|u\in\Lambda\}$, where $u$ is
the leading term of the polynomial $h_u=u-f_u$ in $k\langle X_+
\rangle$.  Let
$$
A=E_k(M,X,(u))=k\langle (\{m_{_j}\}\cup X)_+| S_1 \rangle
$$
where $\{m_{_j}\}_J$ is a basis of $M$,
$S_1=\{u=f_u+(u),u\in\Lambda, \  xm_j=[xm_j], m_jx=[m_jx], \
m_jm_l=0,x\in X,j,l\in J\}$ and $\{(u)|u\in\Lambda\}\subseteq M$.
For convenience, we also call $\{(u)|u\in\Lambda\}$ a factor set of
$B$ in $M$.

We order the set $(\{m_{_j}\}_J\cup X)^+$ also by the deg-lex
order which extends $<_B$ and satisfies $x>m_j$ for any $x\in X$
and $j\in J$.

In $S_1$, the possible compositions are related to the following
ambiguities:
$$
w_1=w, \ m_ju, \ um_j, \ m_jm_lm_n, \ m_jm_lx, \ xm_jm_l, \ xm_jx',
\ m_jxm_l $$ where $w$ is an ambiguity appeared in $R$.

For $w_1=m_jm_lm_n, \ m_jm_lx, \ xm_jm_l, \ xm_jx', \ m_jxm_l$, by
noting that $[[m_jm_l]m_n]=[m_j[m_lm_n]]=0, \
[[m_jm_l]x]=[m_j[m_lx]]=0, \ [[xm_j]m_l]=[x[m_jm_l]]=0, \
[[m_jx]m_l]=[m_j[xm_l]]=0$ and $[[xm_j]x']=[x[m_jx']]$ for any
$j,l,n\in J, \ x,x'\in X$, the corresponding compositions are
trivial modulo $S_1$.

For $w_1=m_ju, \ u=x_1\cdots x_t$, we have
\begin{eqnarray*}
&&(m_jx_1-[m_jx_1],u-(f_u+(u)))_{w_{_1}}=-[m_jx_1]x_2\cdots x_t+m_j(f_u+(u))\\
&\equiv& -[\cdots[[m_jx_1]x_2]\cdots x_t]+m_jf_u+m_j(u) \equiv
-[m_ju]+[m_jf_u]\equiv 0 \ \ mod(S_1,w_1)
\end{eqnarray*}

Similarly, for $w_1=um_j$, we have
$(u-(f_u+(u)),xm_j-[xm_j])_{w_{_1}}\equiv 0 \ \ mod(S_1,w_1)$.

Since $R$ is a minimal Gr\"{o}bner-Shirshov bases, all
compositions in $R$ are only intersection ones.

Now, for $w_1=w=u_1c=du_2, \ u_1,u_2\in \Lambda, \ c,d\in X^+$, we
have
\begin{eqnarray*}
&&(u_1-(f_{u_{_1}}+(u_1)),u_2-(f_{u_{_2}}+(u_2)))_{w_{_1}}
=-(f_{u_{_1}}+(u_1))c+d(f_{u_{_2}}+(u_2))\\
&=& (df_{u_{_2}}-f_{u_{_1}}c)+(d(u_2)-(u_1)c)
\end{eqnarray*}

Since the composition $(h_{u_{_1}},h_{u_{_2}})_{w_{_1}}$ in $R$ is
trivial modulo $R$, it follows that in $k\langle X_+\rangle$,
\begin{eqnarray}\label{e2}
(h_{u_{_1}},h_{u_{_2}})_{w_{_1}}=df_{u_{_2}}-f_{u_{_1}}c=\sum\alpha_i
a_i (u_i'-f_{u_{_i}'}) b_i
\end{eqnarray}
where each $\alpha_i\in k, \ a_i,b_i\in X^{*}, \ u_i'\in\Lambda$
and $\overline{a_i (u_i'-f_{u_{_i}'}) b_i}=\overline{a_i u_i'
b_i}<w_1$. Therefore, in $k\langle (\{m_j\}\cup X)_+\rangle$, we
have
\begin{eqnarray*}
&&(u_1-(f_{u_{_1}}+(u_1)),u_2-(f_{u_{_2}}+(u_2)))_{w_{_1}}=
(df_{u_{_2}}-f_{u_{_1}}c)+(d(u_2)-(u_1)c)\\
&=&\sum\alpha_i a_i (u_i'-f_{u_{_i}'}-(u_i')) b_i+\sum\alpha_i a_i
(u_i') b_i-((u_1)c-d(u_2))
\end{eqnarray*}
where $\overline{a_i (u_i'-f_{u_{_i}'}-(u_i')) b_i}=\overline{a_i
u_i' b_i}< w_1$. From this it follows that
\begin{eqnarray}\label{e3}
(u_1-(f_{u_{_1}}+(u_1)),u_2-(f_{u_{_2}}+(u_2)))_{w_{_1}} \equiv
g_{(u_1,u_2)_{_w}}(u) \ \ mod(S_1,w_1)
\end{eqnarray}
where $g_{(u_1,u_2)_{_w}}(u)=\sum\alpha_i a_i (u_i')
b_i-((u_1)c-d(u_2))\in M$ is a function of $\{(u)|u\in\Lambda\}$.

In fact, by (\ref{e2}), we have an algorithm to find the function
$g_{(u_1,u_2)_{_w}}$.

Therefore, we have the following theorem.

\begin{theorem}\label{st3}
Let $B=k\langle X_+|R \rangle$, where $R=\{u-f_u|u\in\Lambda\}$ is a
 minimal Gr\"{o}bner-Shirshov bases for $B$ and $u$ the leading term
 of the polynomial $h_u=u-f_u$ in
$k\langle X_+ \rangle$. Let $M$ be a $k$-module, $M^2=0$, $M$ a
$B$-bimodule, $\{m_j\}_J$ a $k$-basis of  $M$ and $\{(u)|u\in
\Lambda\}$ a factor set of $B$ in $M$. Let
$A=E_k(M,X,(u))=k\langle ( \{m_{_j}\}_J\cup X)_+| S_1 \rangle $
where $S_1=\{u=f_u+(u) ,xm_j=[xm_j],
m_jx=[m_jx],m_jm_l=0,u\in\Lambda,x\in X,j,l\in J\}$. Then, $S_1$
is a Gr\"{o}bner-Shirshov bases for $A$ if and only if
\begin{eqnarray}\label{e4}
\{g_{(u_1,u_2)_{_w}}(u)|(h_{u_{_1}},h_{u_{_2}})_{w} \mbox{ is a
composition in } R\}=0
\end{eqnarray}
where $g_{(u_1,u_2)_{_w}}(u)$ is defined by (\ref{e3}). Moreover,
if this is the case, $A$  is an extension of $M$ by $B$ in a
natural way.
\end{theorem}
\noindent {\bf Proof.} \ Assume that $S_1$ is a
Gr\"{o}bner-Shirshov bases for $A$. Then, by the previous
statements, for any composition $(U_1,U_2)_{_w}$ in $R$,
$g_{(u_1,u_2)_{_w}}=0$ in $\overline{M}=(M+Id(S))/Id(S)$. Then, by
using Lemma \ref{l2}, the algebra $\overline{M}$ is the same as
$M$. So, the (\ref{e4}) holds. Conversely, if (\ref{e4}) holds,
then it is clear that $S_1$ is a Gr\"{o}bner-Shirshov bases for
$A$.

We need only to prove that $A$  is an extension of $M$ by $B$ in a
natural way if (\ref{e4}) holds. By Lemma \ref{l2}, each element
$r\in A$ can be uniquely written as $r=m+b$, where $m\in M, \ b\in
B$ is $R$-irreducible. Hence, $A=M\oplus B$ as $k$-modules with the
following multiplication: for any $m,m'\in M, \ b,b'\in B$,
$$
(m+b)\cdot (m'+b')=mb'+bm'+(u_{b,b'})+bb'
$$
where $bb'\equiv [bb'] \ mod(R), \ [bb']$ is $R$-irreducible and
$(u_{b,b'})\in M$ is a function of $\{(u)|u\in\Lambda\}$. From this
it follows that $M$ is an ideal of $A$ and $A/M\cong B$ as algebras.
\ \ $\square$

\begin{theorem}\label{st4}
Let $B=k\langle X_+|R \rangle$, where $R=\{u-f_u|u\in\Lambda\}$ is
a minimal Gr\"{o}bner-Shirshov bases for $B$ and $u$ the leading
term of the polynomial $h_u=u-f_u$ in $k\langle X_+ \rangle$. Let
$M, \cal{R}$ be $k$-algebras with $M^2=0$. If $\cal{R}$ is an
extension of $M$ by $B$ and  $\sigma: \ \cal{R}$$/M\rightarrow B,
\ r_x+M\mapsto x$ an algebra isomorphism, then $M$ is a
$B$-bimodule with a natural way: for any $x\in X, \ m\in M$,
\begin{eqnarray*}
x\cdot m=r_xm, \ \ \ m\cdot x=mr_x
\end{eqnarray*}
and there exists a factor set $\{(u)|u\in \Lambda\}$ of $B$ in $M$
such that
$$
\{g_{(u_1,u_2)_{_w}}(u)|(h_{u_{_1}},h_{u_{_2}})_{w} \mbox{ is a
composition in } R\}=0
$$
where $g_{(u_1,u_2)_{_w}}(u)$
 is defined by (\ref{e3}). Moreover, $\cal{R}$
$\cong A$ as algebras, where $A=E_k(M,X,(u))=k\langle
(\{m_{_j}\}_J\cup\{X\})_+| S_1 \rangle $, $\{m_j\}_J$ a $k$-basis of
$M$ and $S_1=\{u=f_u+(u) ,xm_j=[xm_j],
m_jx=[m_jx],m_jm_l=0,u\in\Lambda,x\in X,j,l\in J\}$.
\end{theorem}
\noindent {\bf Proof.} \ It is clear that $M$ is a $B$-bimodule
under the given operations.

For any $u\in\Lambda$, since $r_u+M=r_{f_{_u}}+M$, there exists a
unique $(u)\in M$ such that $r_u=r_{f_{_u}}+(u)$. By noting that
each composition $(h_{u_{_1}},h_{u_{_2}})_{w}$ in $R$,
$(h_{u_{_1}},h_{u_{_2}})_{w}\equiv 0 \ mod(R,w)$, we have that
$g_{(u_1,u_2)_{_w}}(u)=0$ in $M$ by (\ref{e2}).

Now, by using the proof of Theorem \ref{st3}, $S_1$ is a
Gr\"{o}bner-Shirshov bases for the algebra $A$ and then the result
follows. \ \ $\square$

\ \

The following theorem, which gives a complete characterization of an
algebra to be an extension of  $M$ by $B$, follows from Theorem
\ref{st3} and Theorem \ref{st4}.

\begin{theorem}\label{st5}
Let $M,B, \cal{R}$ be $k$-algebras with $M^2=0, \ B=k\langle X_+|R
\rangle$, where $R=\{u-f_u|u\in\Lambda\}$ is a minimal
Gr\"{o}bner-Shirshov bases for $B$ and $u$ the leading term of the
polynomial $h_u=u-f_u$ in $k\langle X_+ \rangle$.  Then $\cal{R}$
is an extension of  $M$ by $B$ if and only if $M$ is a
$B$-bimodule and  there exists a factor set $\{(u)|u\in \Lambda\}$
of $B$ in $M$ such that
$$
\{g_{(u_1,u_2)_{_w}}(u)|(h_{u_{_1}},h_{u_{_2}})_{w} \mbox{ is a
composition in } R\}=0
$$
and $\cal{R}$ $\cong E_k(M,X,(u))=k\langle
(\{m_{_j}\}_J\cup\{X\})_+| S_1 \rangle $, where
$g_{(u_1,u_2)_{_w}}(u)$ is defined by (\ref{e3}), $\{m_j\}_J$ a
$k$-basis of $M$ and $S_1=\{u=f_u+(u) ,xm_j=[xm_j],
m_jx=[m_jx],m_jm_l=0,u\in\Lambda,x\in X,j,l\in J\}$.
\end{theorem}

\section{Applications}

Let $M,B, \cal{R}$ be $k$-algebras with $M^2=0$. The previous
theorems give an answer to how to find the conditions which makes
$\cal{R}$ to be an extension of  $M$ by $B$. We call the condition
(\ref{e3}) the extension condition of $M$ by $B$. In fact, for the
extension condition, it is essential to find the function
$g_{(u_1,u_2)_{_w}}(u)$, which can be obtained by algorithm
(\ref{e2}).

\ \

As results, by using the extension conditions, let us give some
examples. We give the characterization of the extension of $M$ by
$B$ when the $B$ is cyclic algebra, free commutative algebra,
universal envelope of the free Lie algebra and Grassman algebra
 respectively.

\begin{theorem}
Let $M,B, \cal{R}$ be $k$-algebras with $M^2=0$ and $B=k\langle
\{x\}_+|x^n=f(x) \rangle$ a cyclic algebra, where $n$ is a natural
number and $f(x)$ is a polynomial of degree less than $n$ such
that $f(0)=0$. Then, $\cal{R}$ is isomorphic to an extension of
$M$ by $B$ if and only if $M$ is a $B$-bimodule, there exists an
$m\in M$ such that $mx=xm$ and $\cal{R}$ $\cong
E_k(M,x,m)=k\langle (\{m_{_j}\}_J\cup\{x\})_+| S \rangle $, where
$\{m_j\}_J$ is any $k$-basis of $M$ and $S=\{x^n=f(x)+m, \
xm_j=[xm_j], m_jx=[m_jx],m_jm_l=0,j,l\in J\}$.
\end{theorem}
\noindent {\bf Proof.} \ Clearly, $R=\{x^n=f(x)\}$ is  a
Gr\"{o}bner-Shirshov bases for $B$. We need only to consider the
composition in $S$: $(x^n-f(x)-m,x^n-f(x)-m)_w, \ w=x^{n+1}$. Thus,
we obtain the extension condition: $mx=xm$. Now, by Theorem
\ref{st5}, the result follows. \ \ $\square$

\begin{theorem}\label{ex1}
Let $M,B, \cal{R}$ be $k$-algebras with $M^2=0$. Let $X=\{x_i|i\in
I\}, \ I$ a well ordered set and $B=k\langle X_+|x_px_q=x_qx_p, \
p>q, \  p,q\in I \rangle$ the free commutative algebra generated
by $X$. Then, $\cal{R}$ is isomorphic to an extension of  $M$ by
$B$ if and only if $M$ is a $B$-bimodule, there exists a factor
set $\{(x_p,x_q)|p>q, \ p,q\in I \}$ of $B$ in $M$ such that for
any $ p,q,r\in I, \ p>q>r$,
\begin{equation}\label{e31}
(x_q,x_r)x_p-x_p(x_q,x_r)+x_q(x_p,x_r)-(x_p,x_r)x_q+(x_p,x_q)x_r-x_r(x_p,x_q)=0
\end{equation}
and $\cal{R}$ $\cong E_k(M,X,(x_p,x_q))=k\langle  ( \{m_{_j}\}_J\cup
X)_+| S \rangle $, where $\{m_j\}_J$ is any $k$-basis of  $M$ and
$S=\{x_px_q=x_qx_p+(x_p,x_q), \ x_pm_j=[x_pm_j],
m_jx_p=[m_jx_p],m_jm_l=0, p>q, \ p,q\in I ,j,l\in J\}$.
\end{theorem}
\noindent {\bf Proof.} \ Let $R=\{x_px_q=x_qx_p| p>q, \ p,q\in
I\}$. Then, for the deg-lex order on $X^+$, $R$ is clearly a
Gr\"{o}bner-Shirshov bases for $B$ and only one kind of
compositions are in $R$, i.e., $(x_px_q-x_qx_p,x_qx_r-x_rx_q)_w, \
w=x_px_qx_r, \ p,q,r\in I, \ p>q>r$. Then, in $S$, by calculating
the composition
$(x_px_q-x_qx_p-(x_p,x_q),x_qx_r-x_rx_q-(x_q,x_r))_w, \
w=x_px_qx_r, \ p,q,r\in I, \ p>q>r$, we obtain the extension
condition (\ref{e31}). Now, by Theorem \ref{st5}, the result
follows. \ \ $\square$

\ \

The following theorem is a generation of Theorem \ref{ex1}.

\begin{theorem}
Let $M,B, \cal{R}$ be $k$-algebras with $M^2=0$. Let $X=\{x_i|i\in
I\}, \ I$ be a well ordered set, $L=Lie_k( X|[x_px_q]=\sum_{i\in
I}\alpha_{pq}^ix_i, \
 p,q\in I )$ a Lie algebra and  $B=U(L)=k\langle X_+|R \rangle$
 the universal envelope of the Lie algebra
$L$, where each $\alpha_{pq}^i\in k$ and
$R=\{x_px_q=x_qx_p+[x_px_q]| p>q, \ p,q\in I\}$.  Then, $\cal{R}$
is isomorphic to an extension of $M$ by $B$ if and only if $M$ is
a $B$-bimodule, there exists a factor set $\{(x_p,x_q)| p,q\in I
\}$ of $B$ in $M$ such that for any $ p,q,r\in I, \ p>q>r$,
\begin{eqnarray*}\label{e32}
&&(x_q,x_r)x_p+x_q(x_p,x_r)+(x_p,x_q)x_r+([x_qx_r],x_p)+(x_q,[x_px_r])+([x_px_q],x_r) \\
&&-(x_p,x_r)x_q-x_p(x_q,x_r)-x_r(x_p,x_q)=0
\end{eqnarray*}
and $\cal{R}$ $\cong E_k(M,X,(x_p,x_q))=k\langle (\{m_{_j}\}_J\cup
X)_+| S \rangle $, where $(x_p,x_q)=-(x_q,x_p), \ (x_p,x_p)=0$ for
any $p,q\in I$, $\{m_j\}_J$ is any $k$-basis of  $M$ and
$S=\{x_px_q=x_qx_p+[x_px_q]+(x_p,x_q), \ x_pm_j=[x_pm_j],
m_jx_p=[m_jx_p],m_jm_l=0,p>q, \ p,q\in I ,j,l\in J\}$.
\end{theorem}
\noindent {\bf Proof.} \ By \cite{bo06}, for the deg-lex order on
$X^+$, $R$ is a Gr\"{o}bner-Shirshov bases for $B$ and only one
kind of compositions are in $R$, i.e.,
$(x_px_q-x_qx_p-[x_px_q],x_qx_r-x_rx_q-[x_qx_r])_w, \ w=x_px_qx_r,
\ p,q,r\in I, \ p>q>r$. Since, in $S$,
\begin{eqnarray*}
&&(x_px_q-x_qx_p-[x_px_q]-(x_p,x_q),x_qx_r-x_rx_q-[x_qx_r]-(x_q,x_r))_w\\
&=&-(x_qx_px_r+[x_px_q]x_r+(x_p,x_q)x_r)+(x_px_rx_q+x_p[x_qx_r]+x_p(x_q,x_r))\\
&\equiv&-(x_q(x_rx_p+[x_px_r]+(x_p,x_r))+[x_px_q]x_r+(x_p,x_q)x_r)\\
&&+((x_rx_p+[x_px_r]+(x_p,x_r))x_q+x_p[x_qx_r]+x_p(x_q,x_r))\\
&\equiv&-x_qx_rx_p-(x_q[x_px_r]+[x_px_q]x_r)-(x_q(x_p,x_r)+(x_p,x_q)x_r)\\
&&+x_rx_px_q+[x_px_r]x_q+x_p[x_qx_r]+(x_p,x_r)x_q+x_p(x_q,x_r)\\
&\equiv&-(x_rx_q+[x_qx_r]+(x_q,x_r))x_p+x_r(x_px_q+[x_px_q]+(x_p,x_q))\\
&&-(x_q[x_px_r]+[x_px_q]x_r)-(x_q(x_p,x_r)+(x_p,x_q)x_r)\\
&&+[x_px_r]x_q+x_p[x_qx_r]+(x_p,x_r)x_q+x_p(x_q,x_r)\\
&\equiv&-(x_q,x_r)x_p-x_q(x_p,x_r)-(x_p,x_q)x_r+(x_p,x_r)x_q+x_p(x_q,x_r)+x_r(x_p,x_q)\\
&&-[x_qx_r]x_p-x_q[x_px_r]-[x_px_q]x_r+[x_px_r]x_q+x_p[x_qx_r]+x_r[x_px_q]
\end{eqnarray*}
and by using Jacobi identity
$[[x_px_r]x_q]+[x_p[x_qx_r]]+[x_r[x_px_q]]=0$,
\begin{eqnarray*}
&&[x_px_r]x_q+x_p[x_qx_r]+x_r[x_px_q]\\
&\equiv&(x_q[x_px_r]+[[x_px_r]x_q]+([x_px_r],x_q))
+([x_qx_r]x_p+[x_p[x_qx_r]]+(x_p,[x_qx_r]))\\
&&+([x_px_q]x_r+[x_r[x_px_q]]+(x_r,[x_px_q]))\\
&\equiv&[[x_px_r]x_q]+[x_p[x_qx_r]]+[x_r[x_px_q]]
+[x_qx_r]x_p+(x_q[x_px_r]+[x_px_q]x_r\\
&&+(x_p,[x_qx_r])+(x_r,[x_px_q])+([x_px_r],x_q))\\
&\equiv&[x_qx_r]x_p+(x_q[x_px_r]+[x_px_q]x_r
+(x_p,[x_qx_r])+(x_r,[x_px_q])+([x_px_r],x_q))
\end{eqnarray*}
we have the extension condition mentioned in the theorem. Now, by
Theorem \ref{st5}, the result follows. \ \ $\square$

\begin{theorem}
Let $M,B, \cal{R}$ be $k$-algebras with $M^2=0$. Let $X=\{x_i|i\in
I\}, \ I$ be a well ordered set,  $B=k\langle X_+|R \rangle$
 the Grassman  algebra, where  $R=\{x_q^2=0, \ x_px_q=-x_qx_p| p>q, \ p,q\in I\}$. Then,
$\cal{R}$ is isomorphic to an extension of $M$ by $B$ if and only
if $M$ is a $B$-bimodule, there exists a factor set
$\{(x_p,x_q)|p\geq q, \  p,q\in I \}$ of $B$ in $M$ such that for
any $ p,q,r\in I, \ p>q>r$,
\begin{eqnarray*}\label{e33}
&&(x_q,x_q)x_r-x_r(x_q,x_q)+(x_q,x_r)x_q-x_q(x_q,x_r)=0\\
&&(x_r,x_r)x_q-x_q(x_r,x_r)+(x_q,x_r)x_r-x_r(x_q,x_r)=0   \\
&&(x_q,x_r)x_p+(x_p,x_r)x_q+(x_p,x_q)x_r
-x_q(x_p,x_r)-x_p(x_q,x_r)-x_r(x_p,x_q)=0
\end{eqnarray*}
and $\cal{R}$ $\cong E_k(M,X,(x_p,x_q))=k\langle (\{m_{_j}\}_J\cup
X)_+| S \rangle $, where $\{m_j\}_J$ is any $k$-basis of  $M$ and
$S=\{x_q^2=(x_q,x_q), \ x_px_q=-x_qx_p+(x_p,x_q), \ x_qm_j=[x_qm_j],
m_jx_q=[m_jx_q],m_jm_l=0,p>q, \ p,q\in I ,j,l\in J\}$.
\end{theorem}
\noindent {\bf Proof.} \ By \cite{bo06}, for the deg-lex order on
$X^+$, $R$ is a Gr\"{o}bner-Shirshov bases for $B$ and the
possible compositions in $R$ are related to the following
ambiguities: $w_1=x_q^2x_r, \ w_2=x_qx_r^2, \ w_3=x_px_qx_r, \
p,q,r\in I, \ p>q>r$. Then, in $S$, by calculating the
corresponding compositions:
$(x_q^2-(x_q,x_q),x_qx_r+x_rx_q-(x_q,x_r))_{w_1}, \
(x_qx_r+x_rx_q-(x_q,x_r),x_r^2-(x_r,x_r))_{w_2}, \
((x_px_q+x_qx_p-(x_p,x_q),x_qx_r+x_rx_q-(x_q,x_r))_{w_3}$,
respectively, we obtain the extension conditions mentioned in the
theorem. Now, by Theorem \ref{st5}, the result follows. \ \
$\square$

 \ \

\noindent{\bf Acknowledgement}: The author would like to express his
deepest gratitude to Professor L. A. Bokut for his kind guidance,
useful discussions and enthusiastic encouragement when the author
was visiting Sobolev Institute of Mathematics.

\end{document}